\documentclass[11pt]{article}

\usepackage{amsmath,amsfonts,amssymb,graphicx}
\usepackage[margin=3cm]{geometry} 

\newcommand{\pdpd}[2]{\frac{\partial #1}{\partial #2}}
\newcommand{\R}{\mathbb{R}}
\newcommand{\T}{\mathbb{T}}
\newcommand{\Z}{\mathbb{Z}}

\newcommand{\Char}{\operatorname{Char}}
\renewcommand{\Re}{\operatorname{Re}}

\begin{document}

\title{On the generalization of wavelet diagonal preconditioning
to the Helmholtz equation}
\author{Christiaan C. Stolk \\[1ex]
KdV Institute for Mathematics, University of Amsterdam, Amsterdam, \\
The Netherlands, email: C.C.Stolk@uva.nl}
\date{October 2010}

\maketitle
\begin{abstract}
We present a preconditioning method for the multi-dimensional
Helmholtz equation with smoothly varying coefficient.
The method is based on a frame of functions, that approximately 
separates components associated with different singular values
of the operator. For the small singular values, corresponding
to propagating waves, the frame functions are constructed
using ray-theory. A series of 2-D numerical experiments 
demonstrates that the number of iterations required for convergence 
is small and independent of the frequency. In this sense the method 
is optimal.
\end{abstract}

\paragraph{Acknowledgement}
This research was partly funded by the Netherlands Organisation for 
Scientific Research through VIDI grant 639.032.509.

\bigskip

\section{Introduction}

In this paper we will describe a preconditioning method for a class of 
Helmholtz operators
\begin{equation} \label{eq:Helmholtz_op}
  H = -\Delta - \frac{\omega^2}{c(x)^2} + i \frac{\alpha(x) \omega}{L c(x)} ,
\end{equation}
where $c(x)$ and $\alpha(x)$ are positive, smoothly varying coefficents,
$\omega$ is the frequency parameter, and $L$ is a typical length of the 
domain, which is inserted to make $\alpha$ dimensionless.
See
\cite{BrandtLivshits1997,ElmanErnstOLeary2001,ErlanggaOosterleeVuik2006,
Livshits2004,Osei-KuffuorSaad2010,EngquistYing2010PreprintHMat,
EngquistYing2010PreprintMPML}
for examples of preconditioners from the literature.

For elliptic equations there are good preconditioning methods, such as
multigrid and wavelet diagonal preconditioning
\cite{Cohen2003,Dahmen1997}.  Their effect can be explained using
analysis in the position-wave number domain, the so called phase
space. With this point of view we first give a short discussion of
wavelet diagional preconditioning, and then turn to its generalization
to the Helmholtz equation.

Consider, as an example problem an elliptic partial differential
equation $B u = f$, where $B$ is an invertible, symmetric, second
order elliptic operator on a perioc domain (we will not go into the
issues related to boundary conditions.)  The main ingredients of the
method are a basis transformation $F$, given by a wavelet transform,
and an invertible diagonal scaling matrix $W$. The matrix $F^*$ is a
wavelet reconstruction operator, its columns are the elements of the
wavelet basis, denoted by $\psi_\mu$, $F$ is a wavelet decomposition
operator. The elements of the diagonal of $W$ satisfy $w_\mu \sim
2^{-2 |\mu|}$, where $|\mu|$ denotes the scale index of the wavelet, and
by convention $|\mu|=-1$ corresponds to the scaling functions. A
symmetrically preconditioned operator $B$ is then example given by an
expression of the form
\begin{equation} \label{eq:A_precon_basis}
  W^{1/2}  F B F^* W^{1/2}.
\end{equation}
Left- and right preconditioned variants also exist, they are
given by $W F B$ and $B F^* W$. 

The triple $(F^*,W^{-1}, F^*)$ can be viewed as an approximate SVD. A
true SVD $(U,S,V)$ would yield a perfect preconditioner, given by
$S^{-1/2} U' B V S^{-1/2} = I$. Equation (\ref{eq:A_precon_basis})
yields an approximation to the identity in the sense that $W^{1/2} F B
F^* W^{1/2}$ is boundedly invertible, uniformly in parameters such as
the grid constant.

The singular values of partial differential operators (PDO's) with {\em
constant} coefficients follow from Fourier analysis. The symbol of
the operator is defined by
\begin{equation}
  B(\xi) = e^{-i x \cdot \xi}  (B e^{i x \cdot \xi}) .
\end{equation}
The singular functions are complex plane waves and the singular values
$B(\xi_\mu)$, with $\xi_\mu$ the wave vector. 
Wavelets are localized in the Fourier domain around
$\| \xi \| = \pi 2^{|\mu|}$, away from zero and infinity, according to the
vanishing moments and smoothness properties. Here being localized means 
that its Fourier transform decays rapidly away from some bounded region.
Using this it can be shown that they act indeed as approximate singular 
function with value $\sim \pi^2 2^{2 |\mu|}$ \cite{Cohen2003,
DahmenProssdorfSchneider1993}.

For variable coefficient PDO's, or pseudodifferential operators, the
analysis becomes more complicated, but in many cases similar results
hold. The symbol becomes a function of $(x,\xi)$
\cite{AlinhacGerard2007, Hormander1985a,Taylor1981}. If the symbol
doesn't vary too much over the area in phase space where a function is
large, then one can often still show that the operator acts as an
approximate multiplication by the value of its symbol at the phase
space support, or is an approximate singular function. This idea is
well known in semiclassical analysis, see e.g.\ 
\cite{Sjostrand1982,Martinez2002,HelfferSjostrand1986,
  HerauSjostrandStolk2005}, and also occurs in numerical analysis
\cite{TrefethenEmbree2005}. 

With a wavelet basis is associated a tiling of phase space. The phase
space can be divided into disjoint sets $\mathcal{T}_\mu$ that
approximately correspond to the support of the wavelets.
In 1-D, the multi-index $\mu$ has two components, $(j,k)$, and the 
tile is given by $2^{-j}[k,k+1] \times 2^{j} ( [\pi,2\pi] \cup [-2\pi,\pi] )$
for $j\ge 0$ and by $[k,k+1] \times [-\pi,\pi]$ for $j=-1$.
The tiling and the matrix $W$ are {\em adapted} to the operator $B$
in the following sense
\begin{equation} \label{eq:w_factor_symbol_tile_criterion} 
  C^{-1} \le w_{\mu} \sigma_B(x,\xi) \le C  , \qquad 
    \text{for all $(x,\xi) \in \mathcal{T}_\mu$, for all $\mu$.}
\end{equation}
In words: The factor $w_\mu$ compensates for the action of the operator,
which is an approximate multiplication, assuming the symbol doesn't 
vary too much over the tile. 
This explains the choice of the wavelets and the matrix $W$, even if 
it cannot be taken literally and must be supplemented with 
appropriate technical arguments.
The property (\ref{eq:w_factor_symbol_tile_criterion}) will 
serve as the guiding criterion for the design of a preconditioner
for the Helmholtz operator.

Applying the criterion (\ref{eq:w_factor_symbol_tile_criterion}) 
to (\ref{eq:Helmholtz_op}), one finds that the requirement that $F$ is a
{\em basis} transform is very restrictive.  To alleviate this requirement,
methods based on {\em frames} have been proposed \cite{Stevenson2003,
BhowmikStolkPreprint2010,HerrmannMoghaddamStolk2008,HerrmannEtAl2009}.  
If $F$ is a tight frame with $F^* F = I$, one can use a symmetric 
preconditioning of the form
\begin{equation} \label{eq:A_precon_frame}
  F^* W^{1/2} F B F^* W^{1/2} F
\end{equation}
Left- and right-preconditioned variants are given by
\begin{equation}  \label{eq:B_precon_frame_LR}
  B F^* W F , \qquad \text{ and } \qquad
  F^* W F B .
\end{equation}
A second form of left-preconditioning is given by $W F B$. The
latter form leads to an overdetermined system of equations that 
is to be solved in the least-squares sense.
The operators (\ref{eq:A_precon_frame}) and (\ref{eq:B_precon_frame_LR})
are invertible if the diagonal entries of $W$ are real and strictly 
positive. 

The criterion (\ref{eq:w_factor_symbol_tile_criterion}) has important
consequences for the preconditioning of the Helmholtz equation.  Let's
start with the situation that the coefficients $c$ and $\alpha$ are
constant. Clearly a wavelet tiling can't be used in a neighborhood of
the set $\| \xi \| \approx \frac{\omega}{c}$. The wavelet will be of
size $O(\omega)$ in the radial direction, and cover a wide range of
values of the symbol from $|H(\xi) | = O(\omega)$ to $|H(\xi)| =
O(\omega^2)$. Based on the criterion they can only be used some
distance $O(\omega)$ away from this set. Therefore, the generalization
of wavelet diagonal preconditioning to the Helmholtz equation must use
a {\em different} set of basis or frame functions, with a different
phase space tiling.

The tiling must be such that the tiles around $\| \xi \| =
\frac{\omega}{c}$ must be of size at most $O(1)$ in the radial Fourier
direction. It follows that they are of size $O(1)$ in the
corresponding spatial direction. In other words, the frame functions have a
macroscopic spatial extent, also when $\omega$ increases and 
the wavelength becomes
smaller. This conclusion can also be reached from a completely
different way, by considering the flow of information in an iterative
algorithm. The distance the information can travel is proportional to
the size of the support of the frame functions. Since we aim at
convergence in a number of steps independent of $\omega$, and since
the solutions of $H u = f$ contain waves propagating over $O(1)$
distance, there must be frame functions with $O(1)$ size of the
support, and
it is natural that these correspond to the propagating wavenumbers
with $\| \xi \| \approx \frac{\omega}{c}$.  One can also observe that
the number of frame functions with small approximate singular values
is much larger than for the elliptic case, it is on the order of
$O(1)$ time the area of the set $\| \xi \| = \frac{\omega}{c}$. If the
number of points per wavelength is kept constant while $\omega$ may
vary, this leads to $O(N^{d-1})$ basisfunctions with small
eigenvalues, where $N$ is the number of grid points in each direction,
and $d$ is the dimension. For constant coefficients of course a
Fourier basis can be used.

In case of variable coefficients the situation becomes more
complicated.  The tiles with small value of the symbol $| H (x,\xi) |$
must satisfy $\| \xi \| \approx \frac{\omega}{c(x)}$ over an $O(1)$
domain. The wavenumber content of the basisfunction must depend on the
position. The frame functions must therefore be {\em adaptively}
chosen, depending on $c$ and $\omega$.  This is unlike many existing
phase space transforms \cite{Mallat2009}, and also unlike certain
transforms that have been used for hyperbolic PDE (but not for
preconditioning)
\cite{CordobaFefferman1978,Smith1998a,CandesDemanet2005}.

In this paper we show that it is possible to construct such a frame. 
The construction uses the WKB approximation in combination with a
paraxial approximation for the wave equation (because ray theory is
used it can be compared with the work of Brandt and Livshits
\cite{BrandtLivshits1997,Livshits2004} except that they treated the
variable coefficient case only in one dimension.)  Our second main
result concerns the convergence of the preconditioned operator. Using
an implementation in Matlab we find that the number of steps required to
converge becomes bounded by a small number, independent of the
frequency.

While the number of steps is small, the cost per step is still
relatively high. Further research into these transforms might improve
this, as could the fast algorithm of \cite{CandesDemanetYing2007} for
the 3-D case. We find that larger problems can be solved than with the
direct method because the memory requirement is reduced. The
computation is split in two steps, a preparation step and an execution
step, such that the results of the preparation step can be reused for
each right hand side.  We find that the computation time for the
execution step is lightly increased compared to the direct method. The
smoothness requirement on the medium, the inclusion of other boundary
conditions or a damping layer, and the extension to the 3-D case are
topics for further research.

The remainder of the paper is structured as follows.
In section~\ref{sec:Helmholtz_class} some properties of our class
of Helmholtz operators are discussed.
The construction of the preconditioner in the continuous setting
is the topic of 
sections~\ref{sec:construct_1D} and~\ref{sec:construct_nD}.
Some aspects of the discretization and implementation are then explained
in section~\ref{sec:implementation}. 
Section~\ref{sec:num_exp} shows the results of our numerical
experiments. We end with a short discussion.

\section{A class of Helmholtz operators\label{sec:Helmholtz_class}}

The class of Helmholtz operator we study is given in 
(\ref{eq:Helmholtz_op}). We assume that there are constants $C_1$ and
$C_2$ such that
\begin{equation}
  0 < C_1 \le c(x) \le C_2 .
\end{equation}
The imaginary part of $H$ is responsible for damping. 
We assume $\alpha(x)$ is $O(1)$ throughout the domain. So there must 
be $C >0$ such that
\begin{equation}
  \frac{1}{C} \le \alpha(x) \le C .
\end{equation}
We assume the domain is rectangular. For simplicity, the boundary 
conditions are assumed to be periodic.
The spatially varying wave length is given by $\frac{2 \pi c(x)}{\omega}$.
The situation of interest is when the domain is many wave lengths
large, but not so large that a discretization of the Helmholtz equation 
becomes impossible.

With the periodic boundary conditions the part of $H$ without the imaginary 
coefficient, which we denote by
\begin{equation}
  H_1 = - \Delta - \frac{\omega^2}{c(x)^2}
\end{equation}
is real and selfadjoint. It has real eigenvalues that may be close to zero.
The imaginary part influences the singular values
of the operators $H$ and $A$. It leads to the lower bound
\begin{equation}
  \| H u \| \ge \frac{\alpha_{\rm min} \omega}{L c_{\rm max} } ,
\end{equation}
where $\alpha_{\rm min}$ is the minimum of $\alpha$ over the domain,
and $c_{\rm max}$ the maximum of $c$.

We consider a finite difference approximation $A$ of $H$
based on the standard five point stencil for the Laplacian
\begin{equation} \label{eq:discrete_H_op}
  (A u)_{i,j}= \frac{1}{h^2} 
    \left( 4 u_{i,j}-u_{i+1,j}-u_{i-1,j}-u_{i,j+1}-u_{i,j-1} \right)
    + \left( - \frac{\omega^2}{c_{i,j}^2} 
        + \frac{i \alpha_{i,j} \omega}{L c_{i,j}} \right) u_{i,j} .
\end{equation}
The preconditioning method is based on the properties of the
continuous operator, therefore we expect the results to be valid for
other discretizations as well, as long as they approximate the
continuous problem with reasonable accuracy, see also the remarks later on.
Setting $\alpha_{i,j}$ to zero in (\ref{eq:discrete_H_op}) yields the 
finite difference approximation $A_1$ of $H_1$. The matrix $A_1$ is
real and selfadjoint, like $H_1$, so that $A$ satisfies
\begin{equation} \label{eq:A_lower_bound}
  \| A u \| \ge \frac{\alpha_{\rm min} \omega}{L c_{\rm max} }.
\end{equation}
Based on the largest eigenvalue of the discrete Laplacian,
given by $8 h^{-2}$, we find the following estimate for $A$
\begin{equation} \label{eq:A_upper_bound}
  \| A u \| \le \left| \frac{8}{h^2} - \frac{\omega^2}{c_{\rm max}^2} + i
    \frac{\omega}{\alpha_{\rm max} L c_{\rm min} } \right| \, \| u \| .
\end{equation}

The number of points in the domain, and the grid size $h$ are in
practice often chosen based on the number of wavelengths that fit in
the domain. For example using a rule to keep a certain number of
points, say $N_{\rm w}$, per wavelength. Reasonable values
are of the order of 15 or 20. This means
\begin{equation}
  h = \frac{2 \pi c_{\rm min}}{N_{\rm w} \omega} .
\end{equation}
Equations (\ref{eq:A_lower_bound}) and (\ref{eq:A_upper_bound}),
and the assumption that $\alpha$ is not too large, so that \\
$\left| \frac{8}{h^2} - \frac{\omega^2}{c_{\rm max}^2} + i
    \frac{\omega}{\alpha_{\rm max} L c_{\rm min} } \right|
\le \frac{8}{h^2}$
imply the following bound for the condition number
\begin{equation}
  \mathcal{K}(A) < \frac{2 L N_{\rm w}^2 \omega c_{\rm max}}
        {\pi^2 \alpha_{\rm min} c_{\rm min}^2} .
\end{equation}
In other words $\mathcal{K}(A) \sim \omega$ if $\omega \rightarrow \infty$
with the given choice of $h$.

We have made the particular choice that the damping term has coefficient
$O(\omega)$. This leads to the suppression of certain resonance effects, 
where small eigenvalues occurring around specific frequencies lead to 
numerical difficulties. At the same time, the solution to Helmholtz 
equation still consist of propagating waves with propagation distance 
on the order of the domain length.
A damping term $O(\omega^2)$ would lead to a propagation distance on
the order of a fixed number of wavelengths when $\omega \rightarrow
\infty$. The problem would become essentially elliptic.  For elliptic
equations, wavelet or multigrid preconditioning methods lead to a
uniformly bounded operator, i.e.\ the problem of finding a
preconditioner is essentially solved.
The resonance effect that were just mentioned depend on subtle global 
properties of the medium. Some of the small singular values can 
likely be removed by using outgoing radiation
boundary conditions, but it appears that internal wave patterns can
also lead to large values for the norm of the resolvent.
In this respect it could be interesting to study consequences of the 
results from the semiclassical analysis of the Schr\"odinger operator 
and its resonances, see e.g.\ \cite{GerardSjostrand1987}.
It is clear that the ideas of this paper do not take into account the
global properties referred to, so we make no claims concerning the situations
without damping present.

\section{Frame transform adapted to the Helmholtz operator in one dimension%
\label{sec:construct_1D}}

In this section we design a frame transform adapted to the
Helmholtz operator with the domain $[0,1]$ and periodic boundary
conditions. The frame transform will consist of two steps. First we separate
the function into three components that are essentially supported
in three regions of phase space: The large wave numbers,
the positive wave numbers of wavelength scale, and the negative
wave numbers of wave length scale. The small wave numbers
(smaller than wavelength scale) are treated together with the 
wavelength scale wave numbers. For each of these components
a transform corresponding to a tiling is applied:
Fourier tiling for the large wave numbers,
and a {\em modulated Fourier tiling} (associated with
a modulated Fourier transform) for the wave length scale wave numbers.
The modulation parameter is a spatial function depening on the function
$c$, and will cause the required curved space tiling.
The idea of a modulated Fourier transform for preconditioning is new to 
our knowledge.

We first apply a windowing operator
\begin{equation}
  G: u \mapsto \begin{pmatrix} u_1 \\ u_2 \\ u_3 \end{pmatrix}
\end{equation}
using
\begin{equation}
  \hat{u}_1(\xi) = \chi_1(\xi) \hat{u}(\xi) , \qquad
  \hat{u}_2(\xi) = \chi_2(\xi) \hat{u}(\xi) , \qquad
  \hat{u}_3(\xi) = \chi_3(\xi) \hat{u}(\xi)
\end{equation}
The function $\chi_j \in C^\infty(\R)$ to minimize to keep signals
localized in space. They must satisfy $\sum_{j=1}^3 \chi_j^2 =1$, then
this map satisfies $G^* G = I$, so that a tight frame can be constructed.
The function $\chi_1$ will be supported away from the zero set of the 
real part of the symbol, we require $\chi_1(\xi) = 0$ on 
$[-1.5 \omega/c_{\rm min}, 1.5 \omega/c_{\rm min} ]$.
The function $\chi_2$ will be supported around the 
positive zero set of $H_1$.
We require $\chi_2(\xi) = 0$ outside 
$[-0.4 \omega/c_{\rm max}, 2 \omega/c_{\rm min} ]$. Similarly, the function 
$\chi_3$ will be supported around the negative zero set of $H_1$, and 
$\chi_3(\xi) = 0$ outside $[-2 \omega/c_{\rm min}, 0.4 \omega/c_{\rm max} ]$.
Within the overlap regions $[-2 \omega/c_{\rm min}, -1.5 \omega/c_{\rm min} ]$,
$[-0.4 \omega/c_{\rm max}, 0.4 \omega/c_{\rm max} ]$,
$[1.5 \omega/c_{\rm min}, 2 \omega/c_{\rm min} ]$ we use a dilated and 
translated version of a smooth cutoff function $h(x)$, such that
$h(1-x)^2 + h(x)^2 = 1$, and given by \cite{BhowmikStolkPreprint2010}
\begin{equation}
  h(x)  = \left\{
        \begin{array}{ll}
                      0 & \text{if} \quad x \le 0,\\
                      \sin\left(
  \frac{\pi}{2}\frac{e^{-\frac{c}{x}}}{e^{-\frac{c}{x}}+e^{-\frac{c}{1-x}}} 
                      \right)  
                        & \text{if}\quad  0 < x < 1,\\
                      1 &  \text{if} \quad x \ge 1 .
               \end{array} \right.
\end{equation}

\begin{figure}
\begin{center}
\begin{minipage}{6cm}
\begin{center}
(a)\\
\includegraphics[width=6cm]{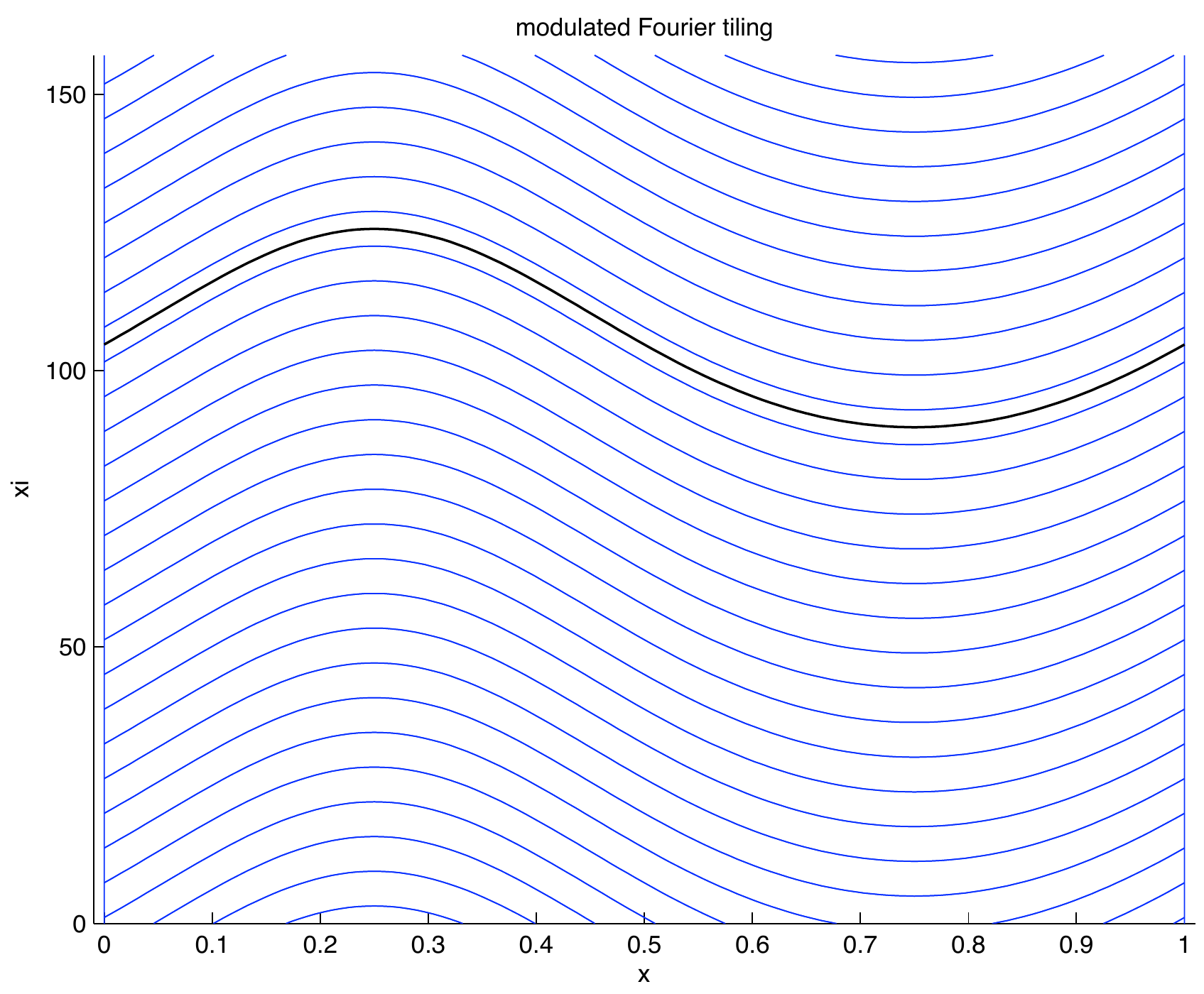}
\end{center}
\end{minipage}\hspace*{8mm}
\begin{minipage}{6cm}
\begin{center}
(b)\\
\includegraphics[width=6cm]{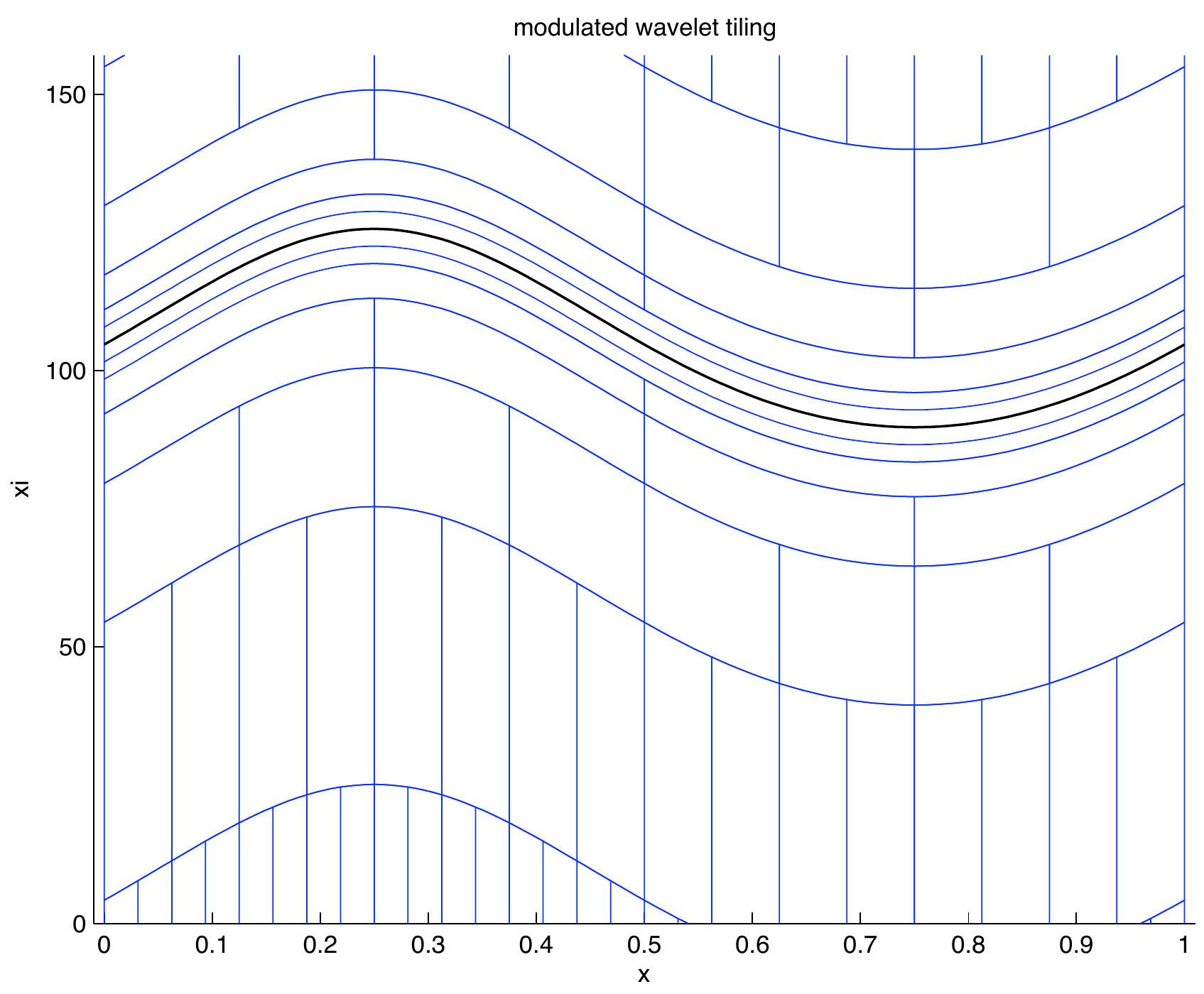}
\end{center}
\end{minipage}
\end{center}
\caption{Example of curved phase space tilings adapted to
$H$ for $\xi \ge 0$. (a) Modulated Fourier tiling; 
(b) modulated wavelet tiling.}
\label{fig:Helmholtz_tiling_1D}
\end{figure}

For each of the $u_j$ we use a basis transform.
For $u_1$ we use the standard Fourier transform
\begin{equation}
  v_j = e^{i x \xi_j}
\end{equation}
The operator $-\partial_x^2 - \omega^2/c(x)^2$ acts on this as
\begin{equation}
  P_1 v_j = (\xi_j^2 - \omega^2/ c(x)^2) v_j
\end{equation}
The preconditioning factor is
\begin{equation}
  \left( ( \xi_j^2 - \omega^2 / c_{\rm mean}^2 )^{2} 
    + \left( \frac{\omega}{c_{\rm mean}} \right)^2 \right)^{-1/2}
\end{equation}
For $u_2$ we use the functions
\begin{equation} \label{eq:modulated_Fourier_1}
  v_j = e^{i \omega T(x) + i x \xi_j }
\end{equation}
with $T(x) = \int_0^x c(s)^{-1} \, ds$ (the travel time), or with a small
linear term added, so that $\omega (T(1) - T(0)) / (2\pi) \in \Z$
and $e^{i \omega T}$ hence satisfies the periodic boundary conditions. The 
values of $\xi_j$ are $\xi_j = 2\pi j$.
The addition of the term $i \omega T(x)$ causes an additive shift in the 
so called instantaneous phase \cite{Mallat2009} by
$\omega \pdpd{T}{x}$. The rectangular tiles associated with the 
Fourier modes $v_j e^{i x \xi_j}$ hence receive an $x$-dependent
shift in the vertical direction of the $(x,\xi)$ plane.
The tiling associated with this basis is given by
\begin{equation}
  0 \le x \le 1 , \qquad
  \frac{\omega}{c(x)} + \xi_j - \pi <
    \xi
  < \frac{\omega}{c(x)} + \xi_j +\pi , \qquad j \in \Z .
\end{equation}
We define the scaling factor for preconditioning by
\begin{equation} \label{eq:preconfac_1D_u2}
  w_j = \left( \left(\frac{\omega}{c(x)} \right)^2
    + \left( 2 \frac{\omega}{c(x)} \xi_j \right)^2 \right)^{-1/2} .
\end{equation}
because we find that the value of the symbol can be uniformly
estimated by $w^{-1}$ on the intersection of the support of 
$u_2$ and the domain of the tile. 
The tiling for an example $c(x)$ is displayed in 
Figure~\ref{fig:Helmholtz_tiling_1D}.
We treat $u_3$ similarly with
\begin{equation} \label{eq:frame_function_1D_u3}
  v_j = e^{-i \omega T(x) + i x \xi_j } ,
\end{equation}
and preconditioning factor given again by (\ref{eq:preconfac_1D_u2}).
Denoting by $F$ the full transform, it follows that $F^* F = I$, i.e.\ 
the transform corresponds to a tight frame.

Instead of (\ref{eq:modulated_Fourier_1}) and
(\ref{eq:frame_function_1D_u3}) it is also possible to use modulated 
wavelets, see the tiling in Figure~\ref{fig:Helmholtz_tiling_1D}.
This can be attractive because a wavelet transform has cost 
$O(N)$, versus $O(N \log N)$ for the Fourier transform.

Some computations using the above methods that were done early on in
this research gave encouraging results.  As the focus of this paper is
on the multi-dimensional case we will not further report on this here.

\section{Frame transform adapted to the higher-dimensional
Helmholtz operator%
\label{sec:construct_nD}}

In this section we devise a frame of functions, assocatiated with
a tiling of phase space adapted to $| H(x,\xi) |$,
in dimension $d = 2$ or higher.
The construction is done in the continuous setting. Discretization and
implementation aspects are discussed in the next section.

Like in the one-dimensional case the most challenging region of phase
space is around the zero-set (characteristic set) of the real part $H_1$, 
given by
\begin{equation}
  \Char(H_1) = \left\{ (x,\xi) \in \T^d \times \R^d \, ; \,
    \| \xi \| = \frac{\omega}{c(x)} \right\} .
\end{equation}
Here $\T^d$ denotes the $d$-dimensional torus, also written as
$[0,1]^d$ with periodic boundary conditions.
In the Fourier space, the direction normal to $\Char(H_1)$, a tile has
size at most $O(1)$, implying that the size in the corresponding spatial 
direction is at least $O(1)$. The tile must follow the curved surface of 
$\Char(H_1)$, a property not present in existing multi-dimensional
transforms, see e.g.\ Mallat \cite{Mallat2009}.

Just like the previous section we will not construct a single, global 
tiling. Instead the first step of the frame transform
is a division into different regions of phase space. Regions associated
with sufficiently large and sufficiently small wave numbers can be 
treated with Fourier preconditioning. The other regions are chosen such
that the part of characteristic set they contain is a graph, say 
$\xi_j = f(x,\xi_1, \ldots, \xi_{j-1}, \xi_{j+1}, \ldots, \xi_n)$.
We call $x_j$ the preferred coordinate. We summarize the remainder of the 
construction. Locally in phase space the Helmholtz equation becomes
equivalent to a first order evolution equation in the preferential direction.
Using theory for hyperbolic pseudodifferential equations an abstract frame
is constructed. To obtain a computable set of functions a WKB 
approximation is made, which finally yields the result.

So the first step in the frame transform is the separation of $u$ into 
components that live in different regions of phase space. 
Localization is done in the Fourier space in radius (scale) and angle,
and in the position space to certain rotated rectangular regions.
There are three wave numbers scales,  
the small wave numbers $\| \xi \| < \omega/c_{\rm max}$, the 
wave length scale wave numbers 
$\omega/c_{\rm max} \le \| \xi \| \le \omega/c_{\rm min}$, and the 
large wave numbers $\| \xi \|  > \omega/c_{\rm min}$.
We will use a map 
$  F_{\rm scale}: u \mapsto \begin{pmatrix} u_1 \\ u_2 \\ u_3 \end{pmatrix}
$
using
\begin{equation} \label{eq:nD_scale_localize_op}
  \hat{u}_1(\xi) = \chi_1(\xi) \hat{u}(\xi) , \qquad
  \hat{u}_2(\xi) = \chi_2(\xi) \hat{u}(\xi) , \qquad
  \hat{u}_3(\xi) = \chi_3(\xi) \hat{u}(\xi)
\end{equation}
with three smooth cutoff functions $\chi_j(\xi)$, $j=1,2,3$ to
separate the small, wavelength scale and large wave numbers. Again 
the cutoff functions will satisfy $\sum \chi_j^2 = 1$ to have 
$F_{\rm scale}^* F_{\rm scale} = I$. By $k_1, k_2$ we denote
constants to separate the low and wave length scale wave numbers. For
$\| \xi \|$ going from $k_1$ to $k_2$ $\chi_1$ goes smoothly from one
to zero and $\chi_2$ smoothly from 0 to 1 according to the standard
cutoff function $h$.  By $k_3, k_4$ we denote constants to separate the
wavelength scale and high wavenumbers.

For further localization in the Fourier space, we will distinguish
several directions in the Fourier space, each direction described by a
unit vector $\alpha \in S^{d-1}$.  Using cutoffs, the support in the
Fourier domain will be restricted to angular wedges around each angle
$\alpha$. This angle will also determine the spatial localization.
We consider new coordinates $(y, z)$, with $y$ in
$\alpha$ direction and $z=(z_1, \ldots, z_{d-1})$ normal to this direction.
We cover the spatial domain with a small number of coordinate boxes 
in these coordinates. The new coordinates will have $y=0$ in the center
of each box.
Because of the periodic boundary conditions, in our 2-D numerical examples
we chose periodic bands. To obtain this, $\alpha$ is set such 
that $\tan(\alpha)$ is a ratio of two small integers. 

After phase space localization and coordinate transform, the preconditioning
can be done according to the following symbol, for which we locally
have $P \sim H_1$
\begin{equation} \label{eq:define_P_0}
  P = 2 \frac{\omega}{c_{\rm mean}}
    \left( \eta - \sqrt{ \frac{\omega^2}{c(y,z)^2} - \zeta^2} \right) .
\end{equation}
Moreover, this can be further modified outside the region of phase space
that the analysis was restricted to. Here this is useful to deal with the 
singularity of the square root at 
$\| \zeta \| = \frac{\omega}{c(y,z)}$. We define a $C^2$ function by
\begin{equation} \label{eq:regularized_one_way}
  S(s_z) = \left\{ \begin{array}{ll} \sqrt{1-s_z^2} 
                      & \text{if $| s_z | \le \sin(\alpha)$} \\
                 \cos(\alpha) - (|s_z| - \sin(\alpha)) \tan(\alpha)
                      & \text{if $| s_z | > \sin(\alpha)$} ,
  \end{array} \right.
\end{equation}
and instead of (\ref{eq:define_P_0}) the following symbol
$P \sim H_1$ will be used
\begin{equation} \label{eq:define_P}
  P = 2 \frac{\omega}{c_{\rm mean}}
   \left( \eta - \frac{\omega}{c(y,z)} S(\frac{\zeta c(y,z)}{\omega}) \right) .
\end{equation}
The function
\begin{equation}
  B(y,z,\zeta,\omega)=\frac{\omega}{c(y,z)} S(\frac{\zeta c(y,z)}{\omega})
\end{equation}
is a symbol (here we ignore that $S$ is only $C^2$ instead of $C^\infty$, 
because a smoother symbol can be designed if the numerical
computations require it). We define
an operator $B(y,z,D_z)$ by
\begin{equation}
  B(y,z,D_z) u(y,z) = 
    (2 \pi)^{-(d-1)} \iint ( B(y,z,\zeta) + \text{l.o.t.} ) 
    e^{i z \cdot \zeta} \hat{u}(y,\zeta) \, d\zeta .
\end{equation}
For later reference we also define a pseudodifferential operator
$B(y,z,D_z,D_t)$ that acts on function of $(z,t)$, with $t$ a time coordinate
\begin{equation}
  B(y,z,D_z,D_t) u(y,z,t) = 
    (2 \pi)^{-d} \iint (B(y,z,\zeta,\omega) + \text{l.o.t.})
    e^{i z \cdot \zeta - i \omega t} \hat{u}(y,\zeta,-\omega) \, d\zeta \, d\omega .
\end{equation}
In these equations $\hat{u}$ defines a Fourier transform w.r.t.\ a subset
of the variables as indicated.
We will assume the lower order terms are such that $B(y,z,D_z)$ and
$B(y,z,D_z,D_t)$ are real and selfadjoint, see e.g.\ \cite{OptRootStolk2010}.
The differential equation assocated with (\ref{eq:define_P}) is 
\begin{equation} \label{eq:one_way}
  \left( \partial_y - i B(y,z,D_z) \right) u(y,z) = 0 .
\end{equation}
It is an evolution equation in $y$.

Equation (\ref{eq:one_way}) is of a well known type, it is called a 
one-way wave equation. The unknown $u$ can be considered as a function 
of $(y,z,\omega)$. After Fourier transform $-\omega \mapsto t$, this results
in a first order hyperbolic pseudodifferential equation results,
see e.g.\ \cite{Hormander1985a,Taylor1981} for the theory of such equations.
One-way wave equations are a tool for the analysis of 
hyperbolic systems of equations, as well as for their numerical computation,
see \cite{OptRootStolk2010} and the references therein.
Because of their nature, they eiter describe waves propagating in the
positive $y$ direction, or in the negative $y$ direction. So called
'turning waves', for which the $y$-component of the propagation direction
changes sign, are not described by such an equation. 
Propagation under large angles with the $y$-axis is also a problem
for numerical methods for one-way wave equations. These usually
become inaccuarate if the angle is too large, e.g.\ larger than 60 
degrees. 

By $E(y,y_0)$ we will denote a solution operator, meaning if 
$u_0 = u_0(z)$ are initial conditions
at $y_0$, then
\begin{equation}
  (\partial_y - i B ) E(y,y_0) u = 0 , \qquad
    E(y_0,y_0) u = u .
\end{equation}
We define the operator $T$, acting on a function of $(y,z)$, by
\begin{equation}
  T u(y,\cdot) = E(y,0) u(y,\cdot)
\end{equation}
From the assumption that the operator $B = B(y,z,D_z)$ is self-adjoint
it follows that this operator is unitary, $T^* T = I$. 
We find that
\begin{equation}
\begin{split}
  (\partial_y - i B) E(y,0) u(y,\cdot)
  = {}& \left( (\partial_y - iB)E(y,0) \right) u(y,\cdot) 
    + E(y,0) \partial_y u 
\\
  = {}& E(y,0) \partial_y u ,
\end{split}
\end{equation}
i.e.\ the following intertwining operator relation holds
\begin{equation}
  (\partial_y - i B) T = T \partial_y .
\end{equation}

Let $\hat{\phi}_\mu(y)$, $\tilde{\phi}_\nu(z)$ be frames, and suppose
the frame $\hat{\phi}_\mu(y)$ is adapted to the operator $\partial_y + 1$
in the sense discussed above, with preconditioning factors $w_\mu$. We 
define a candidate for our frame by
\begin{equation} \label{eq:def_phi_mu_nu}
  \phi_{\mu,\nu}(y,z) = T \hat{\phi}_\mu(y) \tilde{\phi}_\nu(z) .
\end{equation}
Since $T$ is unitary the new frame is also a tight frame with frame 
bound $1$. Indeed,
\begin{equation}
\begin{split}
  \sum_\mu | \langle \tilde{\phi_\mu}, f \rangle | ^2
  = {}& \sum_\mu | \langle T \phi_\mu , f \rangle |^2
  = \sum_\mu | \langle \phi_\mu , T^* f \rangle |^2
\\
  = {}& \langle T^* f, T^* f \rangle = \| f \|^2 .
\end{split}
\end{equation}
We use the intertwining property and the unitarity property, to argue
that the preconditioning weights should be $w_\mu$.
A matrix element satisfies
\begin{equation}
  \langle \phi_{\kappa,\lambda} , (\partial_x - i B) \phi_{\mu,\nu} \rangle
  = \langle \hat{\phi}_\kappa \tilde{\phi}_\lambda ,
    T^* (\partial_y - i B) T \hat{\phi}_\mu \tilde{\phi}_\nu \rangle
  = \langle \hat{\phi}_\kappa , \partial_y \hat{\phi}_\mu \rangle
    \langle \tilde{\phi}_\lambda, \tilde{\phi}_\nu \rangle ,
\end{equation}
hence, taking into account the factor $\frac{2 \omega}{c_{\rm mean}}$ in
(\ref{eq:define_P_0}), $\frac{c_{\rm mean}}{2 \omega}  w_\mu$ is a good 
choice for the preconditioning factors.

Actual application of the frame transform associated with
(\ref{eq:def_phi_mu_nu}) requires the computation of
$E(y,0) \tilde{\phi}_{\nu}(z)$. This amounts to
solving the initial value problem for a first order hyperbolic
pseudodifferential equation. There are various methods to solve
one-way wave equations numerically, see \cite{OptRootStolk2010}
and references. Here we will use a different, semi-analytical
approach, and use the WKB approximation. 
It is assumed that the frame property of the
$\phi_{\mu,\nu}$, and the property that it yields is a good
preconditioner are stable against the (small) perturbations due to the 
errors in the WKB approximation.

For $E(y,0) \tilde{\phi}_\nu(z)$ we therefore look for solutions of the 
form
\begin{equation}
  A(y,z) e^{i \phi(y,z)} .
\end{equation}
The ``initial values'' $\tilde{\phi}_\nu$ must be of the form
$A(z) e^{i \phi(z)}$. Making use of our periodic setting, we choose a 
Fourier basis for $\phi_\nu$ (here $\nu = k \in \Z$)
\begin{equation}
  \tilde{\phi}_k = L_2^{-1/2} e^{i x (2 \pi k/L) } .
\end{equation}
An alternative would be a windowed Fourier frame
$\tilde{\phi}_{k,m} = w(y - k v_0) e^{ i m \zeta_m}$.

Let's start with the equation for $\phi$. 
We write $\phi = \omega T(y,z)$, where $T$ has the dimension of time.
The initial conditions for $T$ are 
$T(0,z) = \frac{2 \pi k }{\omega L } z$.
The WKB method gives the following equation for $T$
\begin{equation} \label{eq:eikonal_oneway_reg}
  \pdpd{T}{y} = S( \frac{c}{\omega} \pdpd{T}{z} ) .
\end{equation}
Without the regularization of (\ref{eq:regularized_one_way}) this becomes
a well-known equation,
\begin{equation} \label{eq:eikonal_oneway}
  \frac{\partial T}{\partial y} 
  - \sqrt{\frac{1}{c^2} - \left( \frac{\partial T}{\partial z} \right)^2 }
  = 0 ,
\end{equation}
the eikonal equation as an evolution equation in the $y$ direction.
Equation (\ref{eq:eikonal_oneway_reg}) can be solved e.g.\ by upwind 
finite differences.

The amplitudes for the equation satisfy a transport equation along
the characteristics \cite{Evans1998}. 
Let $Z(y,z)$ describe a set of characteristics.
The fact that $E$ is unitary provides an additional relation, namely that
\begin{equation} \label{eq:conservation_relation_amplitude}
  A(y,Z(y,z_0)) = A(0,z_0) \left( \pdpd{Z}{z_0} \right)^{-1/2} .
\end{equation}
Instead of computing the characteristics, it is convenient to compute
$Z_0(y,z)$, the inverse of the map $z_0 \mapsto Z(y,z_0)$. This 
function satisfies itself a transport equation that can be solved
along with $T(y,z)$. 

The last step in describing the frame functions is the choice of 
$\hat{\phi}_\mu(y)$. We choose again the Fourier modes, $\mu = j$,
$\hat{\phi}_j(y) = e^{i 2 \pi j /(2L_1)}$.

Having described the frame functions, we look at the resulting 
transform, which, we recall is defined by taking
\begin{equation}
  \langle \phi_{\mu,\nu}, f \rangle
  = \iint \overline{\phi_{\mu,\nu}(y,z)} f(y,z) \, d z \, dy .
\end{equation}
This amounts to two steps. First the map
\begin{equation} \label{eq:cwpt_FIO_step}
  F_1 u(y,\zeta_k) 
  = \int A(y,\zeta_k,z) e^{i\omega T(y,\zeta_k/\omega,y)} u(y,z) \, dz .
\end{equation}
This is a {\em Fourier integral operator}.
Then the map
\begin{equation} \label{eq:cwpt_Fourier_step}
  c_{j,k} = \mathcal{F}_{x \mapsto j} (F_1 u(y,\zeta_k)) .
\end{equation}
We will call the frame transform a Lagrangian wave packet transform
because of the tiling, which, as discussed below, is associated 
with a canonical relation that is by definition a Lagrangian manifold.

We have already discussed that, if $u$ is considered as a function of
$(y,z,\omega)$, and after inverse Fourier transform $-\omega \mapsto
t$, in short, in the time domain, equation (\ref{eq:one_way}) becomes
a first order hyperbolic pseudodifferential equation.
In the time domain the operator $T$ is a Fourier integral operator.
This is useful, because the theory of Fourier integral operators
also includes a geometrical description of the mapping of 
energy in the phase space, which results (for large $\omega$)
in the description of the mapping of singularities. The mapping
of singularities is according to the canonical relation
\cite{Duistermaat1996,Hormander1985b,Treves1980b,Taylor1981}.
The canonical relation for the solution operator to
first order hyperbolic systems is described in terms of the 
bicharacteristics, i.e.\ the characteristics of the eikonal
equation.  We use the notation $B(y,z,\zeta)$ for the 
(principal) symbol of the pseudodifferential
operator $B$.
The bicharacteristics are described by the system of ODE's
\begin{align*}
  \pdpd{z}{y} = {}& - \frac{B}{\zeta} \\
  \pdpd{\zeta}{y} = {}& \frac{B}{z} .
\end{align*}
Based on the unregularized operator (\ref{eq:define_P_0}), the 
explicit equations are
\begin{align*}
  \pdpd{z}{y} = {}& 
    \zeta \left( \frac{\omega^2}{c^2} - \zeta^2 \right)^{-1/2} , \\
  \pdpd{\zeta}{y} = {}& 
    - \frac{\omega^2}{c^3} \pdpd{c}{z} 
      \left( \frac{\omega^2}{c^2} - \zeta^2 \right)^{-1/2} .
\end{align*}
Denote by $(Z(y,z_0,\zeta_0), \Theta(y,z_0,\zeta_0))$ the solution operator
that describes the flow with initial values $(z_0,\zeta_0)$ at $y=0$.
Let $\hat{\mathcal{T}}_\mu$ be the tiling associated with
$\hat{\phi}_\mu$, and $\tilde{\mathcal{T}}_\nu$ be the tiling associated 
with $\tilde{\phi}_\nu$.
This results in the following description of the tiling
\begin{equation}
  \mathcal{T}_{\mu,\nu}
  =
  \{  (y,Z(y,z_0,\zeta_0); \eta+B(y,Z,\Theta), \Theta(y,z_0,\zeta_0) ) \, ; \,
   (z_0,\zeta_0) \in \tilde{\mathcal{T}}_\nu, 
   (y,\eta) \in \hat{\mathcal{T}}_\mu \} .
\end{equation}

\section{Implementation aspects\label{sec:implementation}}

In this section we discuss some aspects of the implementation
of the frame transform and the preconditioner
outlined in the previous section.
The two main steps to be performed for the frame transform are 
the localization in phase space in scale and in angle described
in the first part of this section, and
the transform described in (\ref{eq:cwpt_FIO_step})
and (\ref{eq:cwpt_Fourier_step}). We will comment on the 
discretization of each of those steps, and then say a few 
words on the structure of a program implementing the
preconditioner.

The localization in scale in the Fourier domain was described
in the text around (\ref{eq:nD_scale_localize_op}). This is straightforward
to discretize. We choose $k_1 = 0.2 \omega/c_{\rm max}$,
$k_2 = \omega/c_{\rm max}$, $k_3=1.0 \omega/c_{\rm min}$, 
$k_4 = 1.4 \omega/c_{\rm min}$. The algorithm appears to be insensitive
to the detailed values. 
In connection with the subsampling step described below, 
it is advantageous to use smaller values of $k_4$ if possible.
For the windowing in angle we use $N_{\rm a} = 8$ angles, $-225, -180,
\ldots, 90$ degrees. Three subsequent angles are used to define a
cutoff function in angle, using appropriately translated and dilated
versions of the standard cutoff $h$.

After windowing in the Fourier domain a subsampling step is done.
After localization, the Fourier transform of the signal is supported
in a small subdomain of the original domain. A rectangle is taken 
that contains this subdomain, and the Fourier grid within
the rectangle is used for the further calculations. This results
in a function defined on a substantially coarser grid in space.

The final step associated with the restriction in the phase space
is the rotation to the new spatial coordinates $(y,z)$, and
a resorting of the data into periodic bands. The rotation
amounts to interpolation to a rotated grid. This is done by shear 
rotation which is an exactly invertible operation. It changes the 
grid parameters, but that is not a problem. The shear rotation
is combined with the resorting of the data in periodic bands.
The details of the bands depends on the angle. The total number of 
bands will be denoted by $N_{\rm b}$. The bands are chosen
overlapping, and a cutoff function in the $y$ coordinate
is applied.

For the application of the Fourier integral operator, first 
a preparation step needs to be done computing $T(y,z,p_z)$ and
$A(y,z,p_z)$. The travel time $T(y,z,p_z)$ is computed
using upwind finite differences. These are also used
to compute $Z_0(y,z,p_z)$. $A$ is then computed from 
(\ref{eq:conservation_relation_amplitude}).
The storage for $A(y,\zeta_k,z) e^{i\omega T(y,\zeta_k/\omega,y)}$ is the 
largest storage requirement for the algorithm.
For an $d$ dimensional domain with discretization size $N$ in each
direction, and subsampling toward $\beta N$ in each direction, 
this is about $2 N_{\rm a} ( \beta N)^3$. Some numbers from the 
numerical experiments are given in the next section.

Application of the Fourier integral operator in (\ref{eq:cwpt_FIO_step}) 
is done simply by replacing the integral by a summation.
We have experimented with a partial implementation of the scheme 
described in \cite{CandesDemanetYing2007}, but there were no
significant gains, this was partially related to the large
subsampling we did before computing (\ref{eq:cwpt_FIO_step}).

This leads to the following program structure.
First all preparatory tasks are performed by a routine
{\sc prepare\_helmprec}. The main steps are two call the  
two routines {\sc prepare\_filter\_band} and
{\sc prepare\_lwpt}, that prepare the filtering and banding step,
respectively the Lagrangian wave packet transform.
The results are stored and make it possible to apply the
preconditioner rapidly each time it is called.
The actual execution of the preconditioner is done by a routine
{\sc execute\_helmprec}. A similar subdivision is made, a routine
{\sc execute\_filter\_band} executes the transform
from input data to the filtered and banded data, and the adjoint
of this operation. A routine
{\sc execute\_lwpt} computes the forward version or
the adjoint of the transform described in
(\ref{eq:cwpt_FIO_step}) and (\ref{eq:cwpt_Fourier_step}).
It is straightforward to apply the scaling factors $w_\mu$.

\section{Numerical results\label{sec:num_exp}}

The purpose of our numerical experiments is to study the convergence
of the method. The convergence depends on the problem parameters, the
function $c$, the frequency $\omega$, the damping parameter and the 
right hand side of the equation. Of this, the dependence on $c$
and on $\omega$ are the most interesting, and we will show
some results for different values of $c$ and $\omega$. We have fixed the 
damping to $\alpha = 2 \pi$. A smaller value leads to somewhat slower 
convergence, as one might expect. The right hand side will be chosen as 
a random array. We find that the number of steps to convergence depends
very little on which random right hand side is chosen.

The convergence speed is also influenced by the parameters in the
algorithm. The discretization was chosen to be 16-20 points per
wavelength, where the minimum wavelength present in the domain was
used. The discretization needs to be fine enough to have accurate
dispersion relation, a finer discretization is not important for the
convergence. We used 8 angular directions, and two bands per
angle. The number of bands can increase due to a combination of large
propagation distances (i.e.\ width of the band in $y$ direction) and
strong gradients, in which case the WKB method might no longer be
valid. We did not observe this with the class of media chosen.  After
the filtering in scale and angle, the computations were done on a
coarse grid of about two points per wavelength, which lead to a large
reduction in cost compared to the original grid.  The number of
iterations will be determined by requiring that the error is reduced
with a factor of $10^{-5}$ in the LSQR method using the
right-preconditioned form in (\ref{eq:B_precon_frame_LR}). The
convergence was observed to be linear, which gives an indication for
the other error reduction factors. Other iterative methods were tried,
such as BiCG and BiCGStab, but those were found to perform considerably
worse. This is probably related to the distribution of the eigenvalues
in the domain, eigenvalues were present near the positive and near
the negative real axis.

\begin{figure}
\begin{center}
\includegraphics[width=36mm]{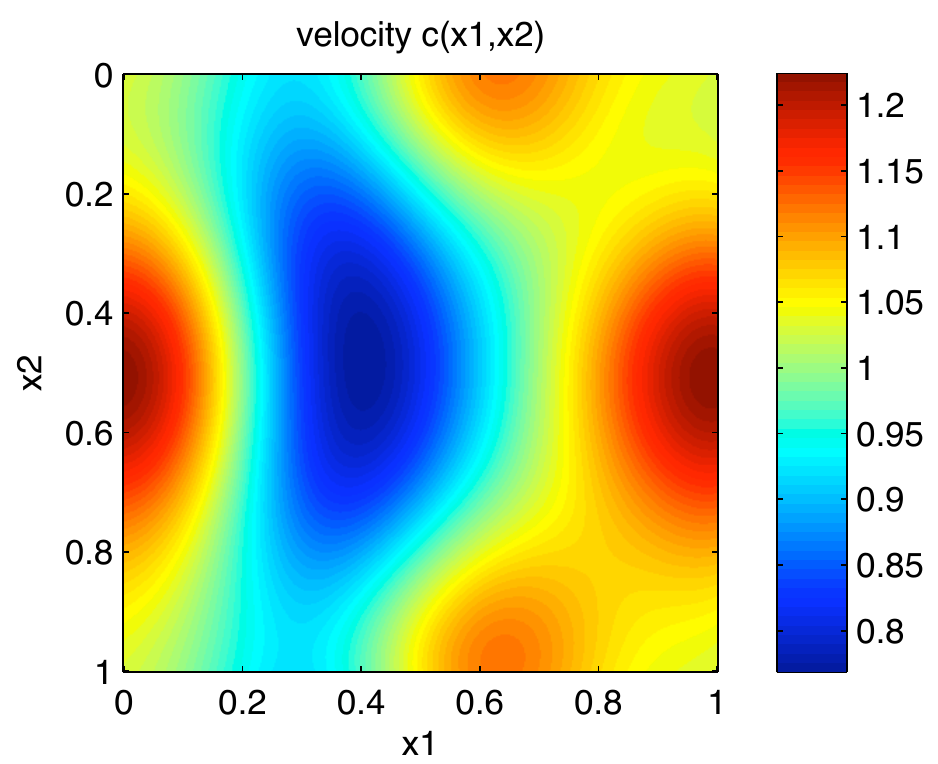}
\includegraphics[width=36mm]{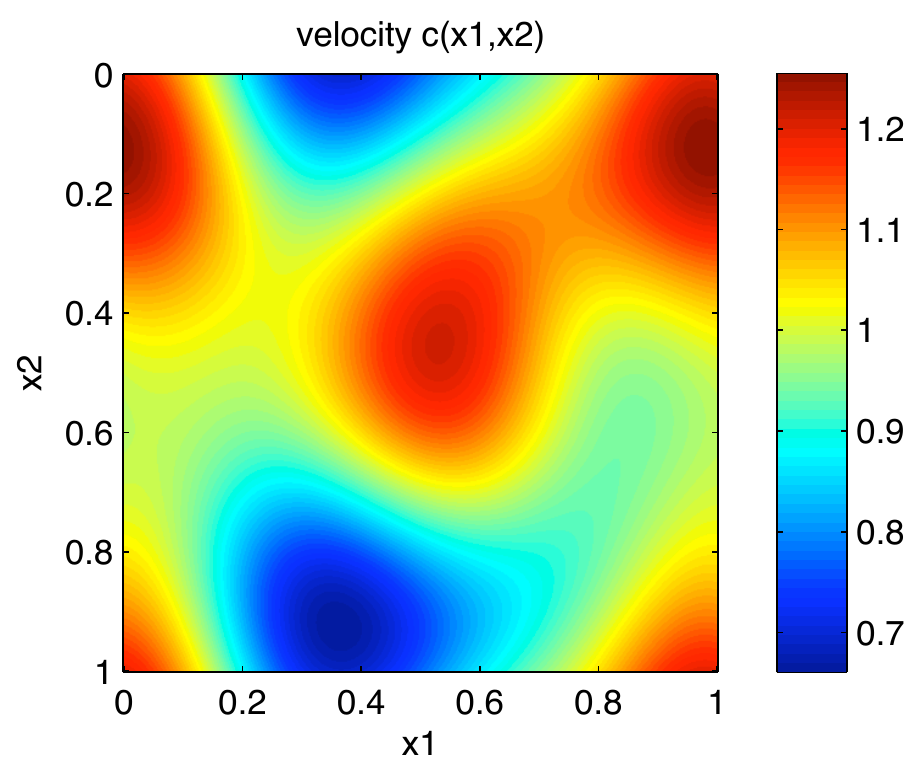}
\includegraphics[width=36mm]{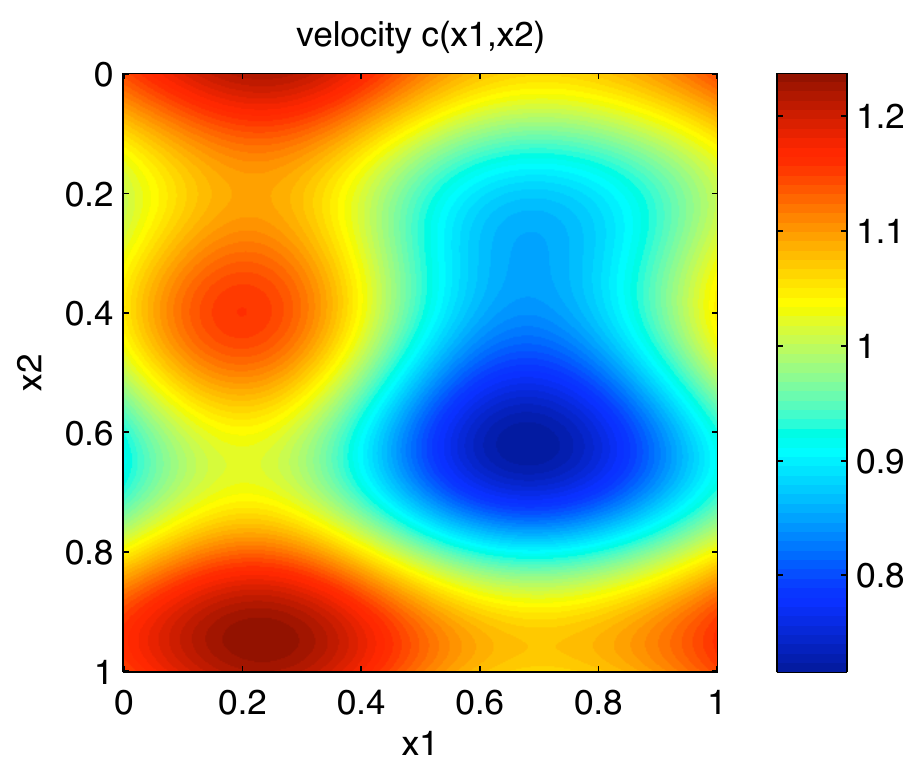}
\includegraphics[width=36mm]{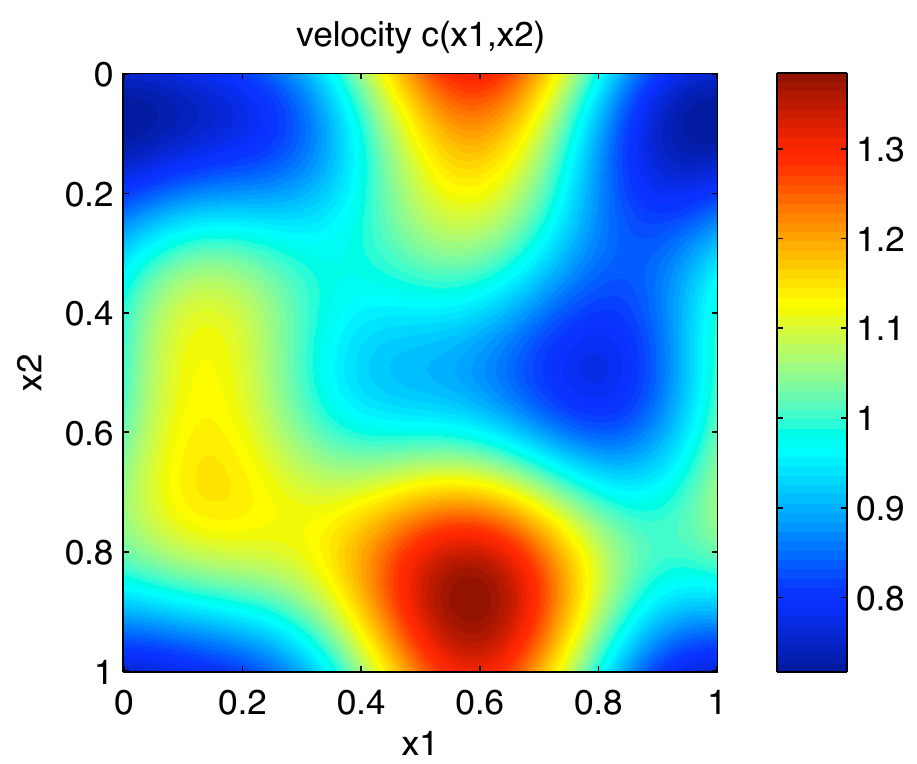}\vspace*{-2ex}
\end{center}
\caption{Some examples for the velocity $c$}
\label{fig:media}
\end{figure}

\begin{figure}[ht]
\begin{center}
\includegraphics[width=52mm]{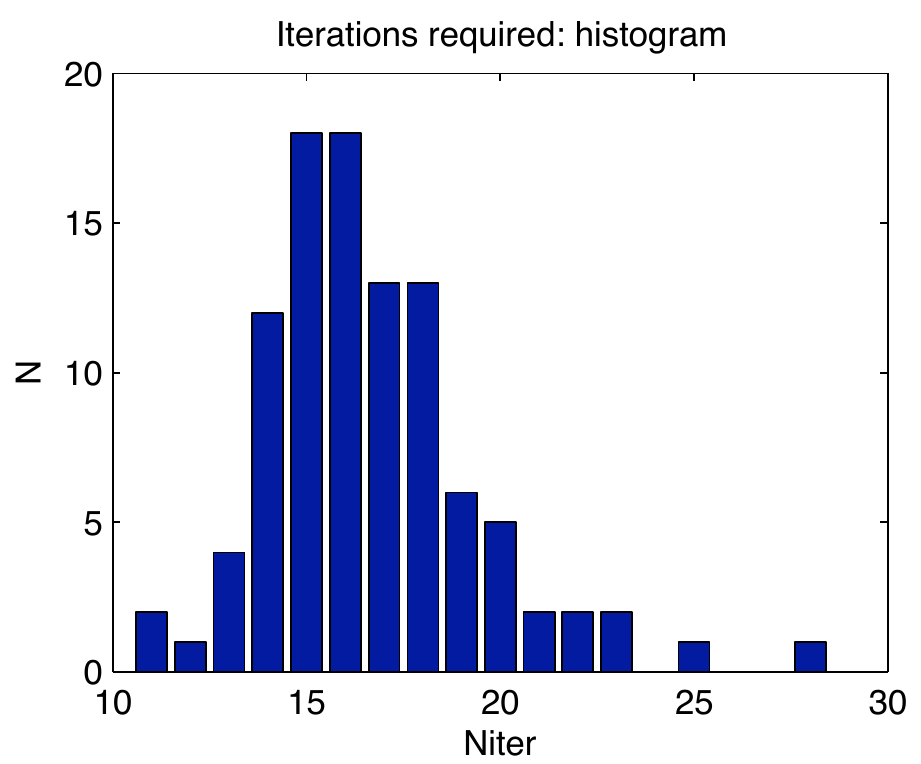}\vspace*{-2ex}
\end{center}
\caption{Histogram describing the number of iterations for convergence
for $\omega = 40\pi$ and 100 different choices of $c$.}
\label{fig:histogram}
\end{figure}

To study the dependence on $c$, we let it vary within a class
of smooth functions, consisting of a sum of a few sines and cosines%
\footnote{%
To be precise, this set is of the form
\begin{equation}
  c(x) = 1 + \Re\left(
      \sum_{j=1}^6 \gamma_{j,1}(a_j+i b_j) 
    e^{2 \pi i (\gamma_{j,2} x_1 + \gamma_{j,3} x_2)} \right)
\end{equation}
where $a_j,b_j$ are random, uniformly distributed in $[-1,1]$,
and the triples $\gamma_j$ take the values
$(0.12,1,0)$, $(0.12,0,1)$, $(0.084,1,1)$, $(0.084,1,-1)$,
$(0.06,2,0)$, $(0.06,0,2)$.%
}, 
The domain will be the unit square. Some examples for $c$ are given in 
Figure~\ref{fig:media}. A set of 100 computations were
performed with different $c$ and $\omega = 40 \pi$, corresponding to
about 20 wavelengths in the domain. A histogram of the results is
given in Figure~\ref{fig:histogram}. There is indeed some spread in
the convergence speed. In general the media with higher velocity
contrasts require more iterations, although we have not found a
precise relation. 

\begin{figure}[ht]
\begin{center}
\includegraphics[width=60mm]{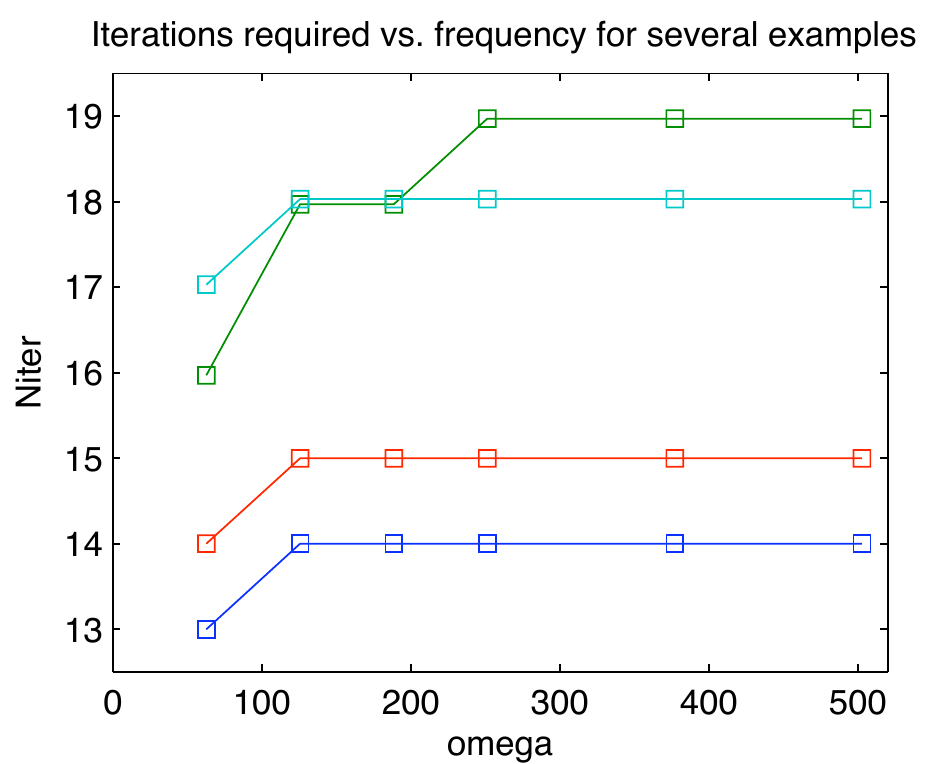}\vspace*{-2ex}
\end{center}
\caption{Number of iterations required vs.\ frequency for four choices
of $c$ within the class described}
\label{fig:iter_vs_freq}
\end{figure}

Our main result concerns the dependence of the convergence on $\omega$.
In Figure~\ref{fig:iter_vs_freq} the number of iterations
as a function of $\omega$ is given for a few choices of $c$.
The main conclusion is that the number of iterations tends to become 
constant for the larger values of $\omega$. For the smaller values of 
$\omega$ the convergence is slightly faster, by one or a few iterations.
The results of Figure~\ref{fig:iter_vs_freq} were quite typical
as $c$ was varied, as can be seen from the comparison with 
Figure~\ref{fig:histogram}


\section{Conclusion and discussion}

We presented a preconditioning method for the multi-dimensional
Helmholtz equation. The preconditioning was based on a frame of functions
designed to have a special phase space tiling, adapted to the Helmholtz 
operator. The tiles follow to a certain degree the level sets of the 
absolute value of the symbol, they are of the general form
\begin{equation}
  \bigcup_{x \in V} \{ x \} \times \Xi(x) .
\end{equation}
Here $\Xi(x)$ is an $x$-dependent subset of the Fourier coordinate space,
$V$ we chosen as a coordinate block. The tiles were curved because of the 
nontrivial dependence of $\Xi(x)$ on $x$.

Our main goal was to establish whether this kind of phase space tiling
could be useful for preconditioning the Helmholtz equation. This is clearly
shown by the numerical experiments. Apart from the convergence,
a second important property was that a large part of the computations could 
be done on very coarse grids.

The generalization to 3-D would of course be of interest.  A question
is whether in 3-D the evaluation of the FIO in
(\ref{eq:cwpt_FIO_step}) can be speeded up by using the results of
\cite{CandesDemanetYing2007}.  It is perhaps useful to point out some
other directions for further development. One issue is the inclusion
of boundary conditions.  For applications like seismic imaging, the
use of a simple absorbing boundary layers can be sufficient. The next
step would be to include a boundary with Dirichlet or Neumann
conditions at a planar or non-planar surface. Less smooth media are
another interesting issue.  We believe the results of this paper form
a strong motivation for such further research.

\bibliographystyle{abbrv} 
\bibliography{refs}

\end{document}